\documentclass[a4paper,12pt]{article}

\usepackage[margin=3cm]{geometry}
\usepackage[pdftex]{graphicx}
\usepackage[utf8]{inputenc}
\usepackage{amsthm}
\usepackage{amsmath}
\usepackage{amssymb}
\usepackage{color}
\usepackage{authblk}
\usepackage{enumerate}
\usepackage{subfigure}

\usepackage[color links=true, backref=page]{hyperref}

\usepackage{soul}

\usepackage{tikz}
\usetikzlibrary{calc}
\usetikzlibrary{backgrounds}
\usetikzlibrary{shapes}
\usetikzlibrary{positioning}
\usetikzlibrary{decorations, decorations.pathmorphing, decorations.text}

\theoremstyle{plain}
\newtheorem{theorem}{Theorem}

\theoremstyle{definition}

\newtheorem{construction}[theorem]{Construction}

\title{The Gray graph is a unit-distance graph}
\author[1]{Leah Wrenn Berman}
\author[2]{G\'{a}bor G\'{e}vay}
\author[3,4,5,6]{Toma\v{z} Pisanski}

\affil[1]{Department of Mathematics \& Statistics,
University of Alaska Fairbanks,
Fairbanks, USA}
\affil[2]{Bolyai Institute, University of Szeged, Szeged, Hungary}
\affil[3]{FAMNIT, University of Primorska, Koper, Slovenia}
\affil[4]{IAM, University of Primorska, Koper, Slovenia}
\affil[5]{Institute of Mathematics, Physics and Mechanics, Ljubljana, Slovenia}
\affil[6]{Faculty of Mathematics and Physics, University of Ljubljana, Ljubljana, Slovenia}

\date{\today}

\begin{document}

\maketitle

\begin{center}
\large{\emph{Dedicated to Dragan Maru\v{s}i\v{c} on the happy occasion\\ of his 70th birthday.}}
\end{center}

\begin{abstract}
In this note we give a construction proving that the Gray graph, which is the smallest cubic semi-symmetric graph, is a unit-distance graph.
\end{abstract}

\noindent \textbf{Keywords:}  
polycirculant,
unit-distance graph,
Gray graph,
ADAM graph,
generalized Petersen graph.

\noindent \textbf{Math.\ Subj.\ Class.\ (2010):} 
05C62,  
05E18,  
20B25,  
05C75,  
05C76,  
05C10.  

\vspace{0.5\baselineskip}


\section{Introduction}

\begin{figure}[!h]
\begin{center}
\includegraphics[width=0.85\textwidth]{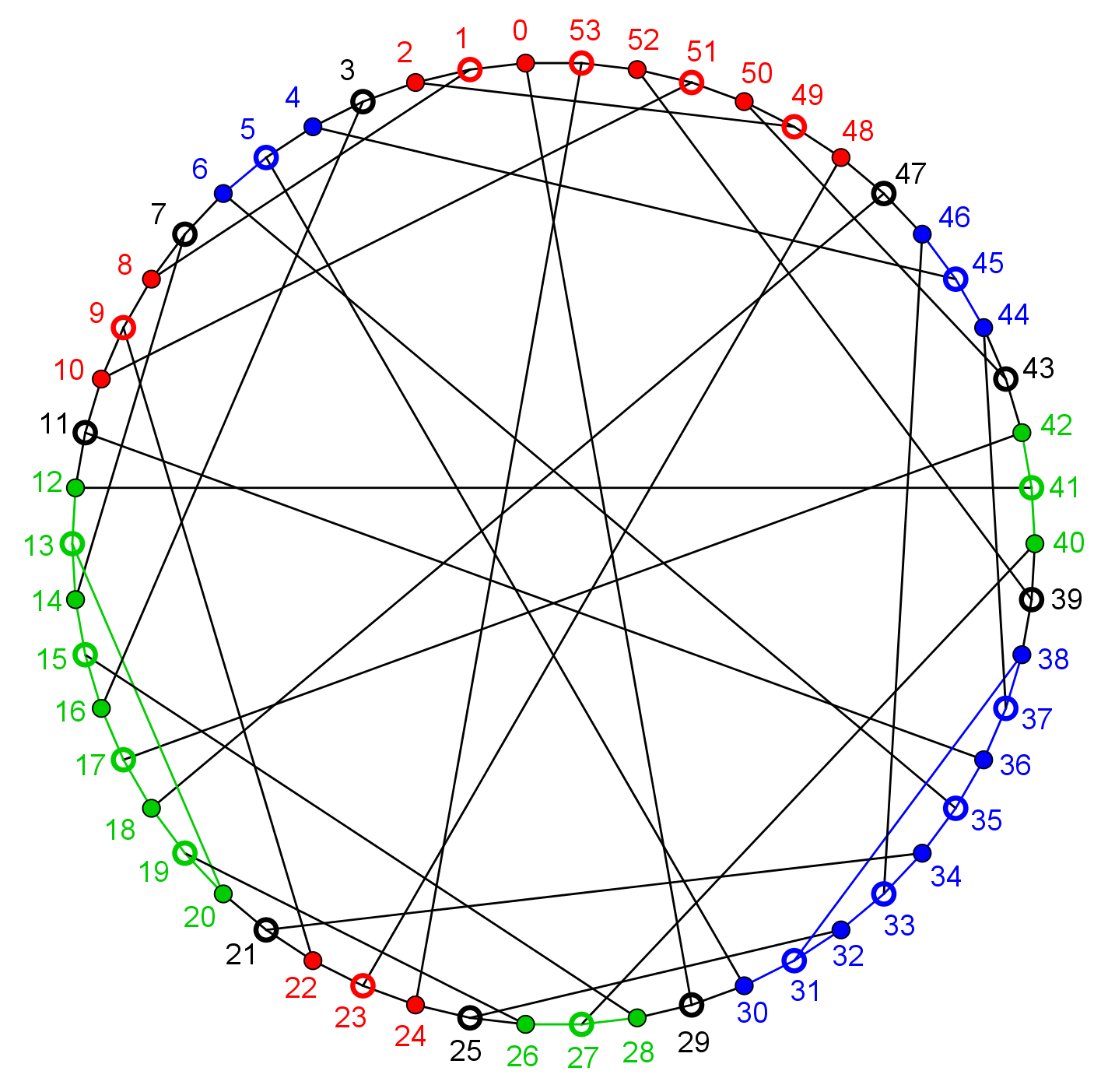}
\caption{The Gray graph with standard vertex labeling, clearly displaying a Hamilton cycle and that it is bipartite, with the bipartition given by hollow and solid vertices. The colour code is explained in the text.}
 \label{fig:HamiltonGrids}
\end{center}
\end{figure}

The well-known \emph{Gray graph}, discovered by Marion C.\ Gray in 1932, has been the subject of 
careful investigation~\cite{MaPi2000, MaPiWi2005, MoPi2007, Pis2007}; see also~\cite[Chapter 6]{PisSer2013}
(for more details about its history, see~\cite{MoPi2007}). It is the smallest trivalent semisymmetric 
graph, meaning that it is the graph with the smallest number of vertices that is 3-regular and edge-transitive but not vertex-transitive~\cite{Bou1972}. It is bipartite and its girth is 8. It is 
Hamiltonian, with a Hamilton cycle displaying a $\mathbb Z_9$ symmetry, as a construction due to Randi\'c shows~\cite{MaPi2000,PiRa2000}. The corresponding LCF notation is $[7, -7, 13, -13, 25, -25]^9$. This can be clearly 
seen in Figure~\ref{fig:HamiltonGrids}.

Clearly, the Hamilton cycle in Figure~\ref{fig:HamiltonGrids} is  symmetric with respect to a cyclic group 
$\mathbb{Z}_3$, shown as 3-fold rotation, generated by
$$
\varrho 
= (0,18,36)(1,19,37)\ldots(17,35,53).
$$

\begin{figure}[!h]
\begin{center}
\includegraphics[width=0.35\textwidth]{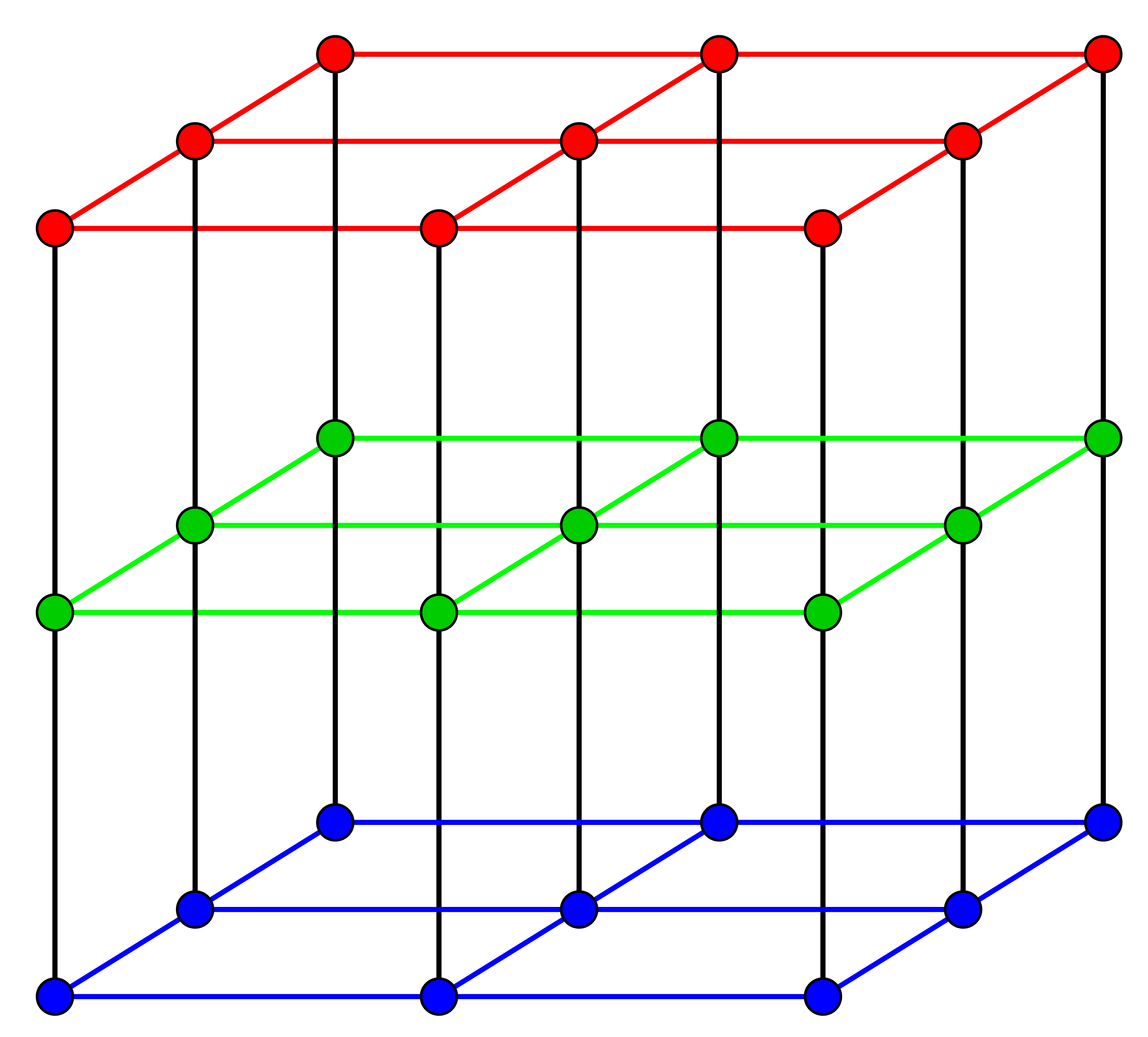}
\caption{The Gray configuration: the standard spatial realization.
             Each coloured horizontal slice represents a $(9_2,6_3)$ subconfiguration. }
 \label{fig:GraySpatial}
\end{center}
\end{figure}

The Gray graph is the Levi graph---that is, the point-line incidence graph---of the \emph{Gray configuration}. This 
is a $(27_{3})$ configuration which occurs as a member of an infinite family of configurations defined by Bouwer 
in 1972~\cite[Section 1]{Bou1972}. Corresponding to this definition, it can be realized as a spatial configuration 
consisting of the 27 points and 27 lines of a $3\times3\times3$ integer grid (see Figure~\ref{fig:GraySpatial}). 
It can be decomposed (in three different ways) into 3 grid-like $(9_2, 6_3)$ subconfigurations lying in 3 parallel 
planes, and a pencil of 9 parallel lines which are perpendicular to these planes. In Figure~\ref{fig:HamiltonGrids} 
we distinguish the vertices representing the points (solid) and lines (hollow) of these subconfigurations by colouring 
red, green and blue, while the lines connecting these copies are represented by black hollow vertices. There is an 
automorphism of the Gray configuration which cyclically permutes these subconfigurations; in terms of the Gray 
graph, the transformation $\varrho$ corresponds just to this automorphism.


\section{A polycirculant unit-distance representations of the Gray graph}


A drawing of a polycirculant graph $\Gamma$ with respect to some semi-regular automorphism $\pi \in \mathbb{Z}_m$ 
is \emph{polycirculant drawing} if $\mathbb{Z}_m$ also acts semi-regularly on the geometric graph, where the edges 
are represented as line segments. Since a cyclic symmetry group in the plane acts as rotations, the vertices of 
the same orbit are equally spaced on a circle. If, in addition, all edge segments are of the same length, the drawing 
is a \emph{polycirculant unit-distance representation} of $\Gamma$ in the Euclidean plane. There are many cases of 
polycirculant drawings in the past, although this is the first time we introduce this term formally. For instance, see 
\cite{ZHP2012} where unit-distance polycyclic drawings of generalized Petersen graphs are presented. Although 
generalized Petersen graphs and $I$-graphs are bicirculants (see e.g.~\cite{BoPiZi2005}), some of them can be 
viewed as tetracirculants, such as the Nauru and ADAM graphs in~\cite{MaPi2018}.

In what follows we give a construction which proves that there exists a polycirculant unit-distance representation of 
the Gray graph with respect to the $\mathbb{Z}_3$ symmetry determined by $\varrho$. It has two degrees of freedom.

\noindent
\begin{construction} \ 

\begin{enumerate}[{\sc{Step }(}1{)}]
\item \label{step:1}
Consider the Levi graph of the $(9_2,6_3)$ subconfiguration of the Gray configuration (see Figure~\ref{fig:grid}a). 
Construct a polycirculant unit-distance representation of it (Figure~\ref{fig:grid}b), e.g.\ in the following way. First 
arrange 6 hollow points so as to form the set of vertices of a regular hexagon. Then choose a radius $r$, and draw 
circles with this radius around each of the 6 points. By suitably pairing these circles, their points of intersection 
yield the desired solid points. Denote this representation by $G_0$. It has threefold rotational (in fact, dihedral) 
symmetry. In addition, it has one degree of freedom, the size of the radius $r$.
\begin{figure}[!h]
\begin{center} 
\subfigure[\hskip -7pt]{
\includegraphics[width=0.275\textwidth]{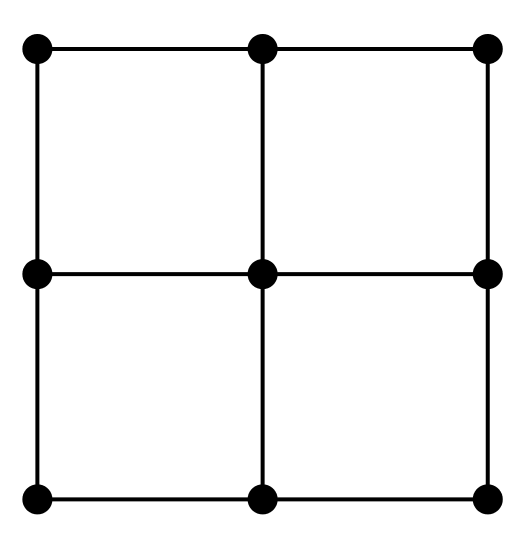}} \hskip 30pt  
\subfigure[\hskip -8pt]{
\includegraphics[width=0.35\textwidth]{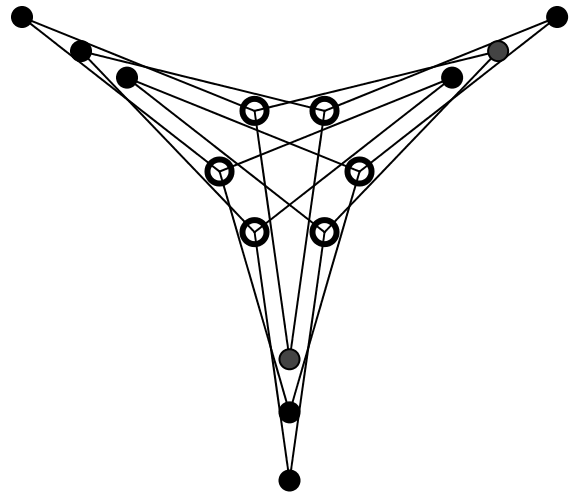}}
\caption{(a) A $(9_2,6_3)$ subconfiguration of the Gray configuration; 
             (b) $G_0$: a polycirculant unit-distance representation of the Levi graph of this subconfiguration. }
\label{fig:grid}
\end{center}
\end{figure}

\item Take a vector star which consists of three unit vectors  that  emanate from a common point whose endpoints 
coincide with the vertices of an equilateral triangle. Denote it by $S$.

\item \label{step:3} Start from $G_0$, and form its three translated images using the vectors of $S$ as translation 
vectors. Denote these translates by $G_b$, $G_g$ and $G_r$. Note that the orientation of the vector star relative 
to $G_0$ determines a new, second degree of freedom of the construction; indeed, the mutual position of the three 
translates of $G_0$ changes when rotating the vector star (see Figure~\ref{fig:Step3}).

\begin{figure}[!h]
\begin{center}
\includegraphics[width=0.85\textwidth]{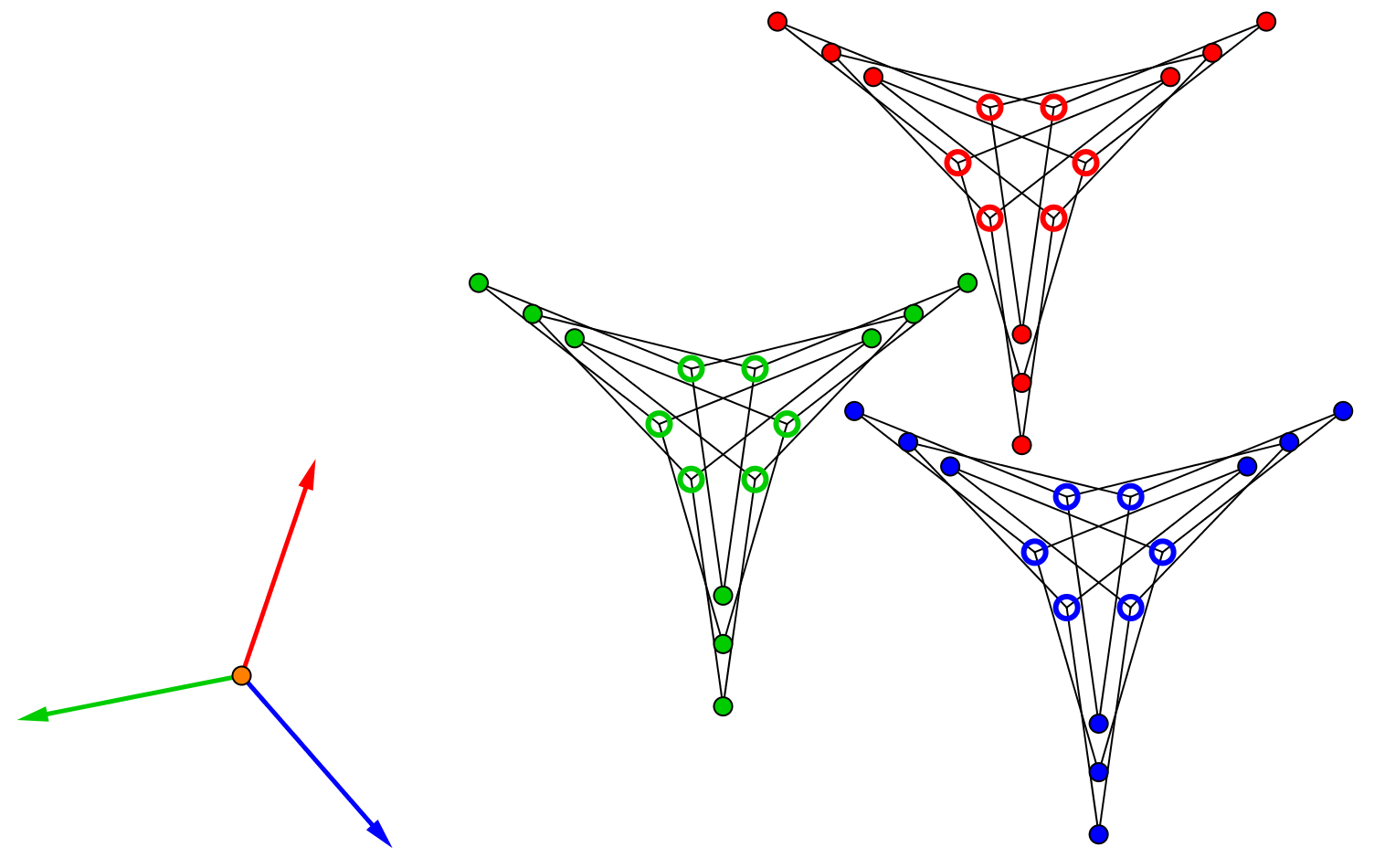}
\caption{Illustration for Step (\ref{step:3}).}
 \label{fig:Step3}
\end{center}
\end{figure}

\item Take the set-union 
\begin{equation*} \label{eq:union}
G_b \cup G_g \cup G_r \cup V_{p}(G_0), 
\end{equation*}
where $V_p(G_0)$ denotes the set of vertices of $G_0$ representing the points of the $(9_2,6_3)$ configuration 
(thus these vertices are denoted by solid vertices). Denote the (disconnected) geometric graph obtained in this way 
by $U$. 

\item In this graph $U$, switch the notation of the 9 points coming from $V_{p}(G_{0})$ 
from solid points to hollow points. Note that, at this stage, $U$ contains $3 \times 9 =27$ solid points and $(3 
\times 6) + 9 = 27$ hollow points, plus 54 edges, coming only from the three subgraphs $G_{b}$, $G_{r}$, and $G_{b}$.

\item Observe that due to Step \eqref{step:3}, in the vicinity of each of the 9 isolated vertices of $U$ there are precisely 
3 solid points at a unit distance, one from each of $G_b$, $G_g$ and $G_r$. Make this graph connected by inserting 
$9 \times 3$ new edges between the isolated vertices and their neighbours of valency 2. The resulting graph now has 
$54 + 27 = 81$ edges and $27$ solid points and $27$ hollow points.
\end{enumerate}
\end{construction}

\begin{figure}[!h]
\begin{center}
\includegraphics[width=0.95\textwidth]{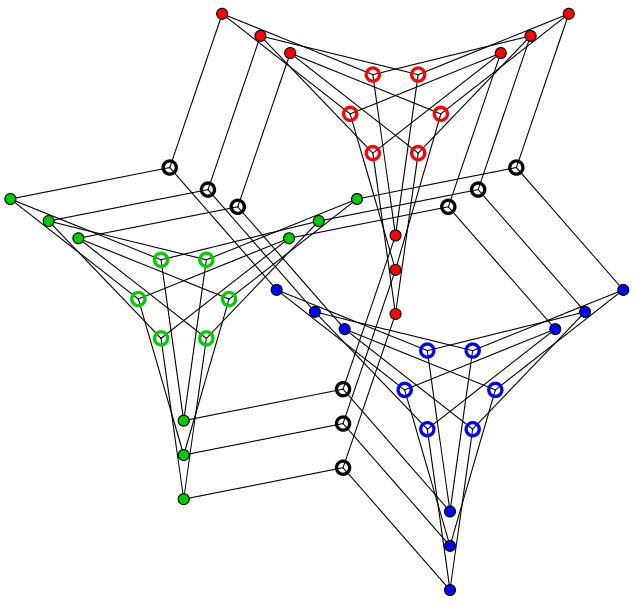}
\caption{A polycirculant unit-distance representation of the Gray graph.}
\label{fig:UD3}
\end{center}
\end{figure}

The final unit-distance graph is shown in Figure \ref{fig:UD3}. It is both polycirculant with $\mathbb{Z}_{3}$ rotational 
symmetry and unit-distance; in addition, it is isomorphic to the Gray graph.
It has two degrees of freedom, as described in Steps~\eqref{step:1} and~\eqref{step:3}, corresponding to a choice 
of unit distance and a choice of angle of rotation of the three coloured subgraphs with respect to the black hollow points.

We note that there is a connection between unit-distance graphs and point-circle isometric configurations
(here ``isometric'' means that all circles in such a configuration are of the same size; see, for instance~\cite{GP14}). 
In our particular case this means that drawing 27 unit circles centred at the hollow vertices of the Figure~\ref{fig:UD3} 
forms, together with the solid vertices, an isometric point-circle realization of the Gray configuration.

Finally, we also note that there are certain similarities between drawing polycyclic configurations and their polycirculant Levi graphs.


 \section*{Acknowledgements}
G\'abor G\'evay is supported by the Hungarian National Research, Development and Innovation Office, OTKA grant 
No.\ SNN 132625.
Toma\v{z} Pisanski is supported in part by the Slovenian Research Agency (research program P1-0294 and research 
projects J1-1690, N1-0140, J1-2481).

\bibliography{GrayPaper}
\bibliographystyle{plain}

\end{document}